\documentclass[11pt]{elsarticle}
\usepackage{amsmath,amssymb,amsthm}
\usepackage{hyperref}
\DeclareFontFamily{U}{mathb}{\hyphenchar\font45}
\DeclareFontShape{U}{mathb}{m}{n}{
	<-6> mathb5 <6-7> mathb6 <7-8> mathb7
	<8-9> mathb8 <9-10> mathb9
	<10-12> mathb10 <12-> mathb12
}{}
\DeclareSymbolFont{mathb}{U}{mathb}{m}{n}
\DeclareMathSymbol{\llcurly}{\mathrel}{mathb}{"CE}
\DeclareMathSymbol{\ggcurly}{\mathrel}{mathb}{"CF}
\DeclareMathOperator{\interior}{int}

\theoremstyle{plain}
\newtheorem{theorem}{Theorem}[section]
\newtheorem{lemma}[theorem]{Lemma}

\newtheorem{proposition}[theorem]{Proposition}
\newtheorem{corollary}[theorem]{Corollary}

\theoremstyle{definition}
\newtheorem{definition}[theorem]{Definition}
\newtheorem{remark}[theorem]{Remark}
\newtheorem{example}[theorem]{Example}

\newcommand{\dda}{\mathord{\mbox{\makebox[0pt][l]{\raisebox{-.4ex}{$\downarrow$}}$\downarrow$}}}
\newcommand{\dua}{\mathord{\mbox{\makebox[0pt][l]{\raisebox{.4ex}{$\uparrow$}}$\uparrow$}}}
\newcommand{\ua}{\mathord{\uparrow}}
\newcommand{\da}{\mathord{\downarrow}}
\newcommand{\rom}[1]{\rm{\uppercase\expandafter{\romannumeral #1}}}

\begin{document}

\begin{frontmatter}




\title{The Directed Probabilistic Powerspace\tnoteref{t1}}
\tnotetext[t1]{Research supported by NSF of China (Nos. 11871353, 12001385).}

\author{Xiaolin Xie}
\ead{xxldannyboy@163.com}

\author{Hui Kou\corref{cor}}
\ead{kouhui@scu.edu.cn}

\author{Zhenchao Lyu}
\ead{zhenchaolyu@scu.edu.cn}
\cortext[cor]{Corresponding author}
\address{Department of Mathematics,\ Sichuan University,\ Chengdu,610064, China }

\begin{abstract}
Probabilistic powerdomain in domain theory plays an important role in modeling the semantics of  nondeterministic  functional programming languages with probabilistic choice.\ In this paper,\ we extend the notion of  powerdomain to directed spaces,\ which is equivalent to the notion of the $T_0$\ monotone-determined space\ \cite{EN2009}.\  We construct the probabilistic powerspace of the directed space,\ which is definded as a free directed space-cone.\ In addition,\ the relationships between our construction  and classical probabilistic powerdomain are studied.\\\



\end{abstract}

\end{frontmatter}

\section{Introduction}
The probabilistic powerdomain was first definded by Saheb-Djahromi [26] in 1980. Its purpose is to provide a mathematical model for the semantics of functional programming languages with probabilistic choice.~Jones proved that if $X$ is a domain, then the probabilistic powerdomain of
$X$ can be realized as the space of valuations, which is still a domain in her PhD thesis [16]. It
has since then been studied extensively by a number of authors, such as, Plotkin[15], Graham[9],
Heckmann [10, 11], Goubault and Jung [7], Tix and Keimel [17, 31], Lyu and Kou [20]. However,
the category of domains together with Scott continuous maps is not cartesian closed. We are
struggling to find a topological representation of the extended probabilistic powerdomain for an
arbitrary dcpo X, but so far there has been no perfect result.
Topological Cones are definded by Keimel [18] in 2006, and his intention is to provide domain
theoretical tools to deal with situations where probabilistic features occur together with ordinary
nondeterminism. In short, over a $T_0$~topological space $X$, he considers the structures that are
close to vector spaces but asymmetric in the sense that elements do not have additive inverses,
accordingly, scalar multiplication is restricted to nonnegative real numbers, besides, both these two operations are jointly continuous.
Directed spaces are introduced by Kou [34] independently, which is equivalent to the $T_0$
monotone determinded space in [4]). A directed space is a $T_0$ space whose topology can be
determined by monotone convergent nets. Especially, every dcpo with the Scott topology is a
directed space, the category of directed spaces with continuous maps is a cartesian closed category
[34]. Hence directed space can be regarded as an extended model of Domain Theory. A natural
question arises: what is a directed space-version extended probabilistic powerspace over a directed
space? To answer this question, we have to do two things. The first one is to define a directed
space-version cone by using Keimel's topological cone, just like the dcpo-cone. The second one is
to construct a free object over any directed space. In this paper we give a concrete construction
of extended probabilistic powerspace over any directed space. The paper is organized as follows:
First to define the notion of directed space-cone and then give its concrete constructure in section
3. In section 4, we give some relations between the extended powerdomain over a dcpo $X$ and
the extended powerspace over $X$.

\section{Preliminaries}
Now,\ we introduce the concepts needed in this article.\ On\ domain\ theory,\ topology,\ and category theory,\ see\ \cite{AJ,  GHK,  Mac71}.\ Let\ $P$\ be a nonempty set.\ A relation\ $\leq$\  on\ $P$\ is called a partial order,\ if\ $\leq$\ satisfies reflexivity\ ($x\leq x$),\ transitivity\ ($x\leq y\ \&\ y\leq z\Rightarrow x\leq z$)\ and antisymmetry\ ($x\leq y\ \&\ y\leq x\Rightarrow x=y$).\ $P$\ is called a partial ordered set(poset) if $\ P\ $is endowed with some partial order\ $\leq$.\ Given\ $A\subseteq P$,\ denote\ $\da A=\{x\in P: \exists a\in A,\  \ x\leq a\}$,\ $\ua A=\{x\in P: \exists a\in A,\ a\leq x\}$.\ We say\ $A$\ is a lower set\ (upper set)\ if\ $A=\da A$ ($A=\ua A$).\ A nonempty set\ $D\subseteq P$\ is called a directed set if each finite nonempty subset of \ $D$\ has upper bound in\ $D$.\ Particularly,\ a poset is called a directed complete poset if each directed subset has a supremum(denoted by\ $\bigvee D$),\ abbreviated as\ dcpo.\ The subset\ $U$\ of poset\ $P$\ is called a\ Scott\ open set if\ $U$\ is an upper set and for each directed set\ $D\subseteq P$,\ which\ $\bigvee D$\ exists and belongs to\ $U$,\ then $U\cap D\not=\emptyset$.\ The set of all Scott open sets of poset\ $P$\ is a toplology on\ $P$,\ which is called the Scott topology and denoted by\  $\sigma(P)$.\ Suppose\ $P,\ E$\ are two posets,\ a fuction\ $f:P\longrightarrow E$\ is called\ Scott\ continuous if it is continuous respect to\ Scott\ topology\ $\sigma(P)$\ and\ $\sigma(E)$.

All topological spaces in this paper are required \ $T_0$\ separation.\ A net of a topological space\ $X$\ is a map\ $\xi: J\longrightarrow X$,\ here\ $J$\ is a directed set.\ Thus,\ each directed subset of a poset can be regarded as a net,\  and its index set is itself.\ Usually,\ we denote a net by\ $(x_j)_{j\in J}$\ or\ $(x_j)$.\ Let\ $x\in X$,\ saying\ $(x_j)$\ converges to\ $x$,\ denote by\ $(x_j)\rightarrow x$\ or\ $x\equiv \lim x_j$,\ if\ $(x_j)$\ is eventually in every open neighborhood of\ $x$,\ that is,\ for each given open neighborhood \ $U$\ of\ $x$,\ there exists\ $j_0\in J$\ such that for every\ $j\in J$,\ $j\geq j_0\Rightarrow x_j\in U$.

Let\ $X$\ be a\ $T_{0}$\ topological space,\ its topology is denoted by\ $O(X)$,\ the specialization order on\ $X$\ is defined as follows:
$$\forall x,\ y \in X,\ x\sqsubseteq y  \Leftrightarrow x\in \overline{\{y\}}$$
here,\ $\overline{\{y\}}$\ means the closure of\ $\{y\}$.
From now on,\ the order of a\ $T_{0}$\ topological space  always indicates the specialization order\ "$\sqsubseteq$".\ Here are some basic properties of specialization order.
\begin{proposition}{\rm\cite{AJ, GHK}}
For a\ $T_{0}$\ topological space\ $X$,\ the following are always true:
	\begin{enumerate}
		\item For each open set\ $U\subseteq X,\ U=\uparrow U;$
		\item For each closed set\ $A\subseteq X,\ A=\downarrow A;$
		\item Suppose\ $Y$\ is another\ $T_{0}$\ topological space,\ and\ $f:X\rightarrow Y$\ is a continuous function from\ $X$\ to\ $Y$.\ Then for each\ $x,\ y\in X$,\ $x\sqsubseteq y \Rightarrow f(x) \sqsubseteq f(y)$,\ that is every continuous fuction is monotone.
	\end{enumerate}
\end{proposition}

Suppose\ $X$\ is a\ $T_{0}$\ space,\ then every directed set\ $D\subseteq X$\ can be regarded as a net of\ $X$,\ we use\ $D\rightarrow x$\ or\ $x\equiv \lim D$\ to represent\ $D$\ converges to\ $x$.\ Define notation
\begin{center}
$D(X)=\{(D,\   x):x\in X,\   D$\ is a directed subset of\ $X$\ and\ $D\rightarrow x \}$.
\end{center}
It is easy to verify that,\ for each\ $x,\   y\in X$,\ $x\sqsubseteq y\Leftrightarrow \{y\}\rightarrow x$.\ Therefore,\ if\ $x\sqsubseteq y$\ then\ $(\{y\},\   x)\in D(X)$.\ Next,\ we give the concept of directed space.

\begin{definition}{\rm\cite{YYK2015}}
	Let\ $X$\ be a\ $T_{0}$\ space.
\begin{enumerate}
\item[(1)] A subset\ $U$\ of\ $X$\ is called a directed open set if
\ $\forall (D,\   x)\in D(X),\ x\in U \Rightarrow D\cap U\neq \emptyset$.\ Denote all directed open sets of \ $X$\ by\ $d(X)$.
\item[(2)] \ $X$\ is called a directed space if each directed open set of\ $X$\ is an open set,\ that is,\ $d(X)=O(X)$.
\end{enumerate}
\end{definition}

\begin{remark}\ \begin{enumerate}
\item[(1)]\ Each open set of a\ $T_{0}$\ space is directed open,\ but the contrary is not necessarily true.\ For example,\ suppose\ $Y$\ is a non-discrete\ $T_1$\ topological space,\ its specialization order is diagonal,\ that is,\ $\forall x,\ y\in Y$,\ $x\sqsubseteq y\Leftrightarrow x=y$.\ Thus,\ all subsets of\ $Y$\ are directed open.\ We notice that\ $Y$\ is non-discrete,\ at least one directed open set is not an open set.
\item[(2)] The definition of directed space here is equivalent to the $T_0$\ monotone determined space defined in\ {\rm \cite{EN2009}}.
\item[(3)] Every poset endowed with the\ Scott\ topology is a directed space\ {\rm\cite{Kou15, YYK2015}},besides,each Alexandroff space is a directed space.\ Thus,\ the directed space extends the concept of the Scott topology.
\end{enumerate}
\end{remark}

Next,\ we introduce the directed continuous function.

\begin{definition}\label{def dc}
	Suppose\ $X,\ Y$\ are two\ $T_{0}$\ spaces.\ A function\ $f:X\longrightarrow Y$\ is called directed continuous if it is monotone and preserves all limits of directed set of\ $X$;\ that is,\ $(D,\   x)\in D(X)\Rightarrow(f(D),\ f(x))\in D(Y)$.
\end{definition}

Here are  some characterizations of the directed continuous functions.
\begin{proposition}\label{prop dc}{\rm\cite{YYK2015}}
	Suppose\ $X,\ Y$\ are two\ $T_{0}$\ spaces.\ $f:X\longrightarrow Y$\ is a function between\ $X$\ and\ $Y$.
	\begin{enumerate}
		\item[(1)]$f$\ is directed continuous if and only if\ $\forall U\in d(Y),\ f^{-1}(U)\in d(X)$.
		\item[(2)]If\ $X,\ Y$\ are directed spaces,\ then \ $f$\ is continuous if and only if it is directed continuous.
	\end{enumerate}
\end{proposition}
	
Now we introduce the product of directed spaces.

Suppose\ $X,\ Y$\ are two directed spaces.\ Let\ $X\times Y$\ represents the cartesian product of\ $X$\ and\ $Y$,\ then we have a natural partial order on it:\ $\forall (x_{1},\ y_{1}),\ (x_{2},\ y_{2})\in X\times Y$,\ $$(x_{1},\ y_{1})\leq(x_{2},\ y_{2})\iff x_{1}\sqsubseteq  x_{2},\ y_{1}\sqsubseteq y_{2}.$$
which is called the pointwise order on\ $X\times Y$.\ Now,\ we define a topological space\ $X\otimes Y$\ as follows:
\begin{enumerate}
	\item The underlying set of\ $X\otimes Y$\ is\ $X\times Y$;
	\item The topology on\ $X\times Y$\ is generated as follows:\ for each given\ $\leq$-\ directed set\ $D\subseteq X\times Y$\ and\ $(x,\ y)\in X\times Y$,\ $$D\rightarrow(x,\ y)\in X\otimes Y\iff\pi_{1}D\rightarrow x\in X,\ \pi_{2}D\rightarrow y\in Y,\   $$
	That is,\ a subset\ $U\subseteq X\times Y$\ is open if and only if for every directed limit defined as above\ $D\rightarrow (x,\ y)$,\ $(x,\ y)\in U\Rightarrow U\cap D\neq\emptyset$.
\end{enumerate}
\begin{theorem}\label{theorem opc}{\rm\cite{YYK2015}}
	Suppose\ $X$\ and\ $Y$\ are two directed spaces.
 \begin{enumerate}
 \item The topological space\ $X\otimes Y$\ defined as above is a directed space and satisfies the following properties:\ the specialization order on\ $X\otimes Y$\ equals to the pointwise order on\ $X\times Y$,\ that is\ $\sqsubseteq =\leq$.

\item Suppose\ $Z$\ is another directed space,\ then\ $f:X\otimes Y\longrightarrow Z$\ is continuous if and only if it is continuous in each variable separately.
\end{enumerate}
\end{theorem}

Denote the category of all directed spaces with continuous functions by\ ${\bf Dtop}$.\ It is proved in\ \cite{YYK2014,YYK2015} that,\ ${\bf Dtop}$\ contains all posets endowed with the Scott topology and\ ${\bf Dtop}$\ is a\ cartesian closed category;\ specifically,\ the categorical products of two directed spaces\ $X$\ and\ $Y$\ are isomorphic to\ $X\otimes Y$.\ So,\ the directed space is an extended framework of\ Domain Theory.

Let\ $P$\ be a\ dcpo,\ and\ $x,\ y\in P$.\ We say\ $x$\ way below $y$,\ if for each given directed set\ $D\subseteq P$,\ $y\leq\bigvee D$\ implies there exists some\ $d\in D$\ such that\ $x\leq d$.\ We write\ $\dda x=\{a\in P:a\ll x\}$,\ $\dua x=\{a\in P:x\ll a\}$.

\begin{definition} A\ dcpo\ $P$\ is called a continuous domain if for each\ $x\in P$,\ $\dda x$\ is directed and\ $x=\bigvee \dda x$.
\end{definition}

\begin{theorem}{\rm\cite{GHK}}
 Suppose\ $P$\ is a continuous\ domain.\ The followings hold:
\begin{enumerate}
\item[(1)] $\forall x,\ y\in P$,\ $x\ll y\Rightarrow \exists z\in P,\ x\ll z\ll y$.
\item[(2)] $\forall x\in P$,\ $\dua x$\ is a\ Scott\ open set.\ Particularly,\ $\{\dua x:x\in P\}$\ is a base of\ $(P,\ \sigma(P))$.
\end{enumerate}
\end{theorem}

C-space was definded by Ern\'{e}\ in\ \cite{CS1991}\ in 1991,\ and it is not hard to verify that each c-space is a directed space.
\begin{definition}\cite{CS1991}
	A\ $T_0$\ topological space $X$\ is a c-space if each\ $x\in X$\ and each open neighborhood\ $U$\ of\ $x$,\ there exists some\ $y\in U$\ such that\ $x\in\interior(\uparrow y)\subseteq U$.
	
\end{definition}

\section{The directed probabilistic powerspaces }
As mentioned above,\ directed space can be regarded as an extended model of\ Domain\ theory,\ just like the work done in article\ \cite{BS2015},\ extending the powerdomain in the category of directed space is very meaningful.\ In this section,\ we will construct the directed probabilistic powerspace of the directed space,\ which is a free algebra generated by the addition and scalar multiplication of the directed space.

\begin{definition}{\rm\cite{GHK}}
	Let\ $X$\ be a topological space,\ a map\ $\mu:\mathcal{O}(X)\longrightarrow[0,\ +\infty]$\ is called a continuous valuation if the followings hold:
	\begin{enumerate}[(1)]
		\item  strictness:\ $\mu(\emptyset)=0$;
		\item monotonicity:\ $V\subseteq U$\ implies\ $\mu(V)\leq\mu(U)$;
		\item modular law:\ $\mu(U)+\mu(V)=\mu(U\cup V)+\mu(U\cap V)$;
		\item continuity:\ for each directed family\ $\mathcal{D}\subseteq\mathcal{O}(X)$,\ $\mu(\sup\mathcal{D})=\sup\{\mu(U):U\in\mathcal{D}\})$.
	\end{enumerate}
	
\end{definition}

\begin{definition}\cite{GHK}
	Let\ $X$\ be a topological space, for each\ $x\in X$,\ we define the point valuation\ $\eta_x:\mathcal{O}(X)\longrightarrow [0,+\infty]$\ from the lattice of open sets of\ $X$,\ $\eta_x(U)=1$\ if\ $x\in U$, and\ $\eta_x(U)=0$\ if\ $x\notin U$.\ A finite linear sum\ $\xi=\sum_{b\in B}r_b\eta_b$\ with\ $0< r_b < +\infty$,\ and defined by\ $\xi(U)=\sum_{b\in U}r_b$\ is called a simple valuation,\ and the set\ $B$\ is called its support.
\end{definition}

Topological cones was defined by Keimel in\ \cite{KEI2006}\ in 2006,\ here we replace the\ $T_0$\ space by directed space to get the definition of directed space-cone without changing any other conditions.
	
\begin{definition}\label{def dcone}
		Let\ $X$\ be a directed space.\ A directed space-cone is a directed space\ $X$\  equipped with a distinguished element\ $0\in X$,\ an addition\ $+:X\otimes X\rightarrow X$,\ and a scalar multiplication\ $\cdot:\mathbb{R}^+\otimes X\rightarrow X$($\mathbb{R}^+$\ endowed with the Scott topology in its usual partial order)\ such that both operations are jointly continuous and the followings are satisfied.
\begin{enumerate}
\item $x+y=y+x$,\ $\forall x,\ y\in X$,
\item $(x+y)+z=x+(y+z)$,\ $\forall x,\ y,\ z\in X$,
\item $0+x=x,\ \forall x\in X$,
\item $(kl)\cdot x=k\cdot(l\cdot x),\ \forall k,\ l\in\mathbb{R}^+,\forall x\in X$,
\item $(k+l)\cdot x=(k\cdot x)+(l\cdot x)$,\ $\forall k,\ l\in\mathbb{R}^+,\ \forall x\in X$,
\item $k\cdot(x+y)=(k\cdot x)+(k\cdot y),\ \forall k\in\mathbb{R^+},\ \forall x,\ y\in X$,
\item $1\cdot x=x,\ \forall x\in X$,
\item $k\cdot 0=0$,\ $\forall k\in\mathbb{R}^+$.
\end{enumerate}

	\end{definition}

\begin{definition}
Suppose\ $(X,\ +,\ \cdot),\ (Y,\ \uplus,\ \ast)$\ are two directed space-cones,\ $f:(X,\ +,\ \cdot)\rightarrow (Y,\ \uplus,\ \ast)$\ is called a directed space-cone homomorphism between\ $X$\ and\ $Y$,\ if\ $f$\ is continuous and\ $f(x+ y)=f(x)\uplus f(y),\ f(a\cdot x)=a\ast f(x)$\ holds for\ $\forall x,\ y\in X,\ a\in\mathbb{R}^+$\ .
\end{definition}

Denote the category of all directed space-cone and directed space-cone homomorphisms by\ ${\bf Dcone}$.\ Then\ ${\bf Dcone}$\ is a subcategory of \ ${\bf Dtop}$.

Next,\ we give the definition of directed probabilistic powerspace.

\begin{definition}\label{def dpp}
  Suppose\ $X$\ is a directed space.\ A directed space\ $Z$\ is called the directed probabilistic powerspace of\ $X$\ if  the following two conditions are satisfied:
\begin{enumerate}
\item[(1)] $Z$\ is a directed space-cone,\ that is the addition\ $+ $\ and the scalar multiplication\ $\cdot$\ exist and they are continuous,
\item[(2)] There is a continuous function\ $i:X\longrightarrow Z$\ satisfies:\ for arbitrary directed space-cone\ $(Y,\ \uplus,\ \ast)$\ and continuous function\ $f:X\longrightarrow Y$,\ there exists a unique directed space-cone homomorphism\ $\bar{f}:(Z,\ +,\ \cdot)\rightarrow (Y,\ \uplus,\ \ast)$\ such that\ $f=\bar{f}\circ i$.
\end{enumerate}
\end{definition}

By the definition above,\ if directed space-cone\ $(Z_1,\ +,\ \cdot)$\ and\ $(Z_2,\ \uplus,\ \ast)$\ are both the directed probabilistic powerspaces of\ $X$,\ then there exists a topological homeomorphism which is also a directed space-cone homomorphism\ $g:Z_1\rightarrow Z_2$.\ Therefore,\ by the means of algebraic isomorphism and topological homeomorphism,\ the directed probabilistic powerspaces of a directed space is unique.\ Particularly,\ we denote the directed probabilistic powerspace of each directed space\ $X$\ by\ $P_P(X)$.

\vskip 3mm

Now,\ we will prove the existence of the directed probabilistic powerspace of each directed space\ $X$\ by way of concrete construction.

\vskip 3mm

	Let\ $X$\ be a directed space.\ Let\ $\mathcal{SV}(X)$\ denotes the set of all simple valuations in\ $X$,\ we have the pointwise order
	\ $\leq$\ on\ $\mathcal{SV}(X):$$$\xi\leq\eta\iff \xi(U)\leq\eta(U),\ \forall U\in\mathcal{O}(X).$$

	Let\ $\mathcal{D}=\{\xi_i\}_{i\in I}\subseteq \mathcal{SV}(X)$\ be a directed set,\ $\xi=\sum_{i=1}^n r_{b_i}\eta_{b_i}\in \mathcal{SV}(X)$\ with support\ $B$.\ We say\ $\mathcal{D}\Rightarrow_P \xi$\ if the following two conditions hold:
	\begin{enumerate}[(1)]
		\item for each\ $b_i\in B$,\ there exists a directed set\ $D_i\subseteq X$\ such that in\ $X$\ with\ $D_i\longrightarrow b_i,i=1,\dots,n.$;
		\item $\forall\ (d_1,\dots,d_n)\in\prod_{i=1}^{n}D_i,\ \forall\  r_{b_i}^{'}<r_{b_i},i=1,\dots,n$,\ there exists some\ $\xi^{'}\in\mathcal{D}$,\ such that\ $\sum_{i=1}^{n} r_{b_i}^{'}\eta_{d_i}\leq\xi^{'}$.
	\end{enumerate}

A subset\ $\mathcal{U}\subseteq \mathcal{SV}(X)$\ is called a\ $\Rightarrow_P$\ convergence open set of\ $\mathcal{SV}(X)$\ if and only if for each directed subset\ $\mathcal{D}$\ of\ $\mathcal{SV}(X)$\ and\ $\xi\in \mathcal{SV}(X)$,\ \ $\mathcal{D}\Rightarrow_{P}\xi\in \mathcal{U}$\ implies\ $\mathcal{D}\cap \mathcal{U}\neq \emptyset$.\ Denote all\ $\Rightarrow_P$\ convergence open sets of\ $\mathcal{SV}(X)$\ by\ $O_{\Rightarrow_{P}}(\mathcal{SV}(X))$\ .

\begin{proposition}\label{prop CX}
	Suppose\ $X$\ is a directed space,\ the following are true:
	\begin{enumerate}
		\item[(1)] $(\mathcal{SV}(X),\ O_{\Rightarrow_P}(\mathcal{SV}(X)))$\ is a topological space,\ abbreviated as\ $CX$.
		\item[(2)] The specialization order\ $\sqsubseteq $\ on\ $CX$\ equals to\ $\leq$.
		\item[(3)] $C(X)$\ is a directed space,\ that is\ $O_{\Rightarrow_P}(\mathcal{SV}(X))=d(CX)$.
		
	\end{enumerate}
\end{proposition}
\noindent{\bf Proof}
	 (1)  Obviously we have\ $\emptyset,\ CX\in O_{\Rightarrow_P}(\mathcal{SV}(X))$.\ If\ $\mathcal{U}\in O_{\Rightarrow_U}(UX)$,\ and\ $\xi\leq\eta,\eta\in\mathcal{U},\ \xi=\sum_{i=1}^{n}r_{b_{i}}\eta_{b_i}$.\ Then it is evident that\ $\{\eta\}\Rightarrow_{P}\xi$,\ since we only need to take\ $D_{i}=\{b_{i}\},\ i=1,\dots,n$.\ Then,\ $\{\eta\}\cap \mathcal{U} \neq
	\emptyset$,\ this means\ $\eta\in \mathcal{U}$,\ and\ $\mathcal{U}$\ is an upper set respect to order\ $\leq_P$,\ $\mathcal{U}=\uparrow_{\leq}\mathcal{U}$.
	
	Let\ $\mathcal{U}_{1},\mathcal{U}_{2}\in O_{\Rightarrow_P}(\mathcal{SV}(X))$,\ and a directed set \ $\mathcal D\subseteq CX$\ with\ $\mathcal{D}\Rightarrow_{P} \xi\in \mathcal{U}_{1}\cap\mathcal{U}_{2}$,\ then,\ there exists\ $\xi_{1}\in \mathcal{D}\cap \mathcal{U}_{1}$\ and\ $
	\xi_{2}\in \mathcal{D}\cap \mathcal{U}_{2}$,\ but\ $\mathcal{D}$\ is directed,\ we have\ $\xi_{3}\in \mathcal{D},\ \xi_{3}\geq\xi_{1},\ \xi_{2}$.\ Then,\ $\xi_{3}\in \mathcal{D}\cap \mathcal{U}_{1}\cap\mathcal{U}_{2}$.\ By the same way,\ we can prove that \ $O_{\Rightarrow_U}(\mathcal{SV}(X))$\ is closed under arbitrary union.\ It follows that\ $O_{\Rightarrow_U}(\mathcal{SV}(X))$\ is a topology.

\vskip 3mm
\
(2) By the proof of\ (1),\ each\ $\Rightarrow_P$\ convergence open set is an upper set respect to\ $\leq$,\ then\ $\downarrow_{\leq}\eta\subseteq\downarrow_{\sqsubseteq}\eta,\forall\eta\in CX$.\ We only need to prove
that\ $\downarrow_{\leq}\eta$\ is a closed set in\ $CX$,\ since\ $\downarrow_{\sqsubseteq}\eta$\ is the minimal closed set containing\ $\eta$.\ Equivantly,\ $CX\setminus\downarrow_{\leq}\eta$\ is a\ $\Rightarrow_{P}$\ convergence open set.

Set\ $\mathcal{U}=CX\setminus\downarrow_{\leq}\eta$.\ Suppose we have\ $\mathcal{D}\Rightarrow_{P}\xi=\sum_{i=1}^{n}r_{b_{i}}\eta_{b_i}\in\mathcal{U}$.\ By contradiction,\ suppose\  $\mathcal{U}\cap\mathcal{D}=\emptyset$,\ that is\ $\forall\xi^{'}\in\mathcal{D},\xi^{'}\leq\eta$.\ However,\ by the definition of\ $\mathcal{D}\Rightarrow_{P}\xi$,\ we have directed sets\ $D_1,\dots,D_n\subseteq X$\ with\ $D_i\rightarrow b_i,i=1,\dots,n.$\ and\ $\forall\ (d_1,\dots,d_n)\in\prod_{i=1}^{n}D_i,\forall r_{b_{i}}^{'}<r_{b_{i}}$,\ there exists some\ $\xi^{'}\in\mathcal{D}$\ such that\ $\sum_{i=1}^{n} r_{b_i}^{'}\eta_{d_i}\leq\xi^{'}$.\ Now,\ we claim that\ $\xi\leq\eta$,\ which contradicts with\ $\xi\in\mathcal{U}$.\
By the definition of pointwise order\ $\leq$,\ $\forall U\in\mathcal{O}(X)$,\ we may assume that\ $b_1,\dots,b_k\in U,0\leq k\leq n$.\  Since\ $D_i\rightarrow b_i,i=1,\dots,k$,\ there exists\ $(d_1,\dots,d_k)\in\prod_{i=1}^{k}D_i$\ such that\ $d_i\in U$.\ For each\ $r_{b_{i}}^{'}<r_{b_{i}},i=1,\dots,k$,\ there exists\ $\xi^{'}\in\mathcal{D}$\ such that\ $\sum_{i=1}^{k} r_{b_i}^{'}\eta_{d_i}\leq\xi^{'}$.\ Note that\ $k\leq n$,\ then,\ $$(\sum_{i=1}^{k} r_{b_i}^{'}\eta_{b_i})(U)=\sum_{i=1}^{k} r_{b_i}^{'}\leq\xi^{'}(U)\leq\eta(U).$$\ But the supremum of the left hand side in the upper inequality is\ $\sum_{i=1}^{k} r_{b_i}=(\sum_{i=1}^{k}r_{b_{i}}\eta_{b_i})(U)=\xi(U)$.

(3) For an arbitrary topological space\ $X$,\ $O(X)\subseteq d(X)$\ holds,\ then\ $O_{\Rightarrow_P}(\mathcal{SV}(X))\subseteq d(CX)$.\ On the other hand,\ according to the definition of\ $\Rightarrow_P$\ convergence topology,\ if directed set\ $\mathcal D\subseteq CX$\ with\ $
	          \mathcal{D}\Rightarrow_{P}\xi$,\ then\ $\mathcal{D}$\ convergents to\ $\xi$ respect to\ $O_{\Rightarrow_P}(\mathcal{SV}(X))$.\ Thus,\ by the definition of directed open set,\ $\mathcal{D}\Rightarrow_{P}\xi \in\mathcal{U}\in d(CX)$\ will imply\ $\mathcal{U}\cap \mathcal{D}\neq \emptyset$.\ Then,\ $\mathcal{U}\in O_{\Rightarrow_P}(\mathcal{SV}(X))$,\ it follows that\ $O_{\Rightarrow_P}(\mathcal{SV}(X))= d(CX)$,\ that is,\ $CX$\ is a directed space.\ $\Box$

\begin{proposition}\label{prop CXdc}
	Suppose\ $X,\ Y$\ are two directed spaces.\ Then function\ $f:CX\rightarrow Y$\ is continuous if and only if for each directed set\ ${\mathcal D}\subseteq CX$\ and\ $\xi\in CX$,\ ${\mathcal D}\Rightarrow_P\xi$\ implies\ $f(\mathcal{D})\rightarrow f(\xi)$.
\end{proposition}

\noindent{\bf Proof} Since\ $\Rightarrow_P$\ convergence will lead to\ $O_{\Rightarrow_P}(\mathcal{SV}(X))$\ topological convergence,\ the necessity is obvious.\ We are going to prove the sufficiency. Firstly,\ we check
that\ $f$\ is monotone.\ If\ $\xi,\ \eta\in CX$\ and\ $\xi\leq\eta$,\ then\ $\{\eta\}\Rightarrow_P\xi$,\ by the hypothesis,\ $\{f(\eta)\}\rightarrow f(\xi)$,\ thus\ $f(\xi)\sqsubseteq f(\eta)$.\ Suppose\ $U$\ is an open set of\ $Y$\ and the directed set\ ${\mathcal D}\Rightarrow_P\xi\in f^{-1}(U)$,\ then\ $f(\mathcal{D})$\ is a directed set of\ $Y$\ and\ $f(\mathcal{D})\rightarrow f(\ua F)\in U$,\ thus\ $\exists\  \xi^{'}\in {\mathcal D}$\ such that\ $f(\xi^{'})\in U$.\ That is,\ $\xi^{'}\in {\mathcal D}\cap f^{-1}(U)$.\ According to the definition of\ $\Rightarrow_P$\ convergence open set,\ $f^{-1}(U)\in O_{\Rightarrow_P}(\mathcal{SV}(X))$,\ that means\ $f$\ is continuous.\ $\Box$\\

Define an addition\ $+$\ on\ $CX:\forall \xi,\ \eta\in CX,\ (\xi +\eta)(U)=\xi(U)+\eta(U),\ \forall U\in\mathcal{O}(X)$.\ Define a scalar multiplication\ $\cdot$\ on\ $CX$\ :$\forall a\in\mathbb{R}^{+},\ \forall \xi\in CX,\ (a\cdot\xi)(U)=a\xi(U)$.\ Next,\ we shall check that these two operations are both continuous.\ Thus,\ $(CX,\ +,\ \cdot)$\ is a directed space-cone.

\begin{theorem} Let\ $X$\ be a directed space.\ Then\ $(CX,\ +,\ \cdot)$\ is a directed space-cone.
\end{theorem}

\noindent{\bf Proof}
By Proposition\ \ref{prop CX},\ $CX$\ is a directed space,\ and by the definition of\ $+$\ and\ $\cdot$,\ these two operations are both monotone.\ According to\ Theorem\ \ref{theorem opc}\ and\ Proposition\ \ref{prop CXdc},\ to prove the continuity of\  $+$,\ we only need to check that for arbitrary fixed\ $\eta=\sum_{j=1}^{m}s_{c_j}\eta_{c_j}\in CX$,\ and a directed set\  $\mathcal{D}=\{\xi_i\}_{i\in I}\subseteq CX$\ with\ $\mathcal{D}\Rightarrow_{P}\xi=\sum_{i=1}^{n}r_{b_i}\eta_{b_i}$, then\  $\mathcal{D}+\eta =\{\xi_i +\eta\}_{i\in I}\Rightarrow_{P}\xi +\eta$.\ By hypothesis,\ $\mathcal{D}\Rightarrow_{P}\xi$,\ there exist\  $D_1,\dots,D_n$,\ with\ $D_i\rightarrow b_i$.\ Taking finite single point sets\ $\{c_j\},\ j=1,\dots,m$.\ For arbitrary\  $(d_1,\dots,d_n,c_1,\dots,c_m)\in\prod_{i=1}^{n}D_i\times\prod_{j=1}^{m}\{c_j\}$,\ and\ $\forall r_{b_i}^{'}<r_{b_i},\ s_{c_j}^{'}<s_{c_j},\ i=1,\dots,n,j=1,\dots,m.$\ Again,\ by the definition of\ $\mathcal{D}\Rightarrow_{P}\xi$,\ there exists some\  $\xi^{'}\in\mathcal{D}$\ such that\ $\sum_{i=1}^{n}r_{b_i}^{'}\eta_{d_{i}}\leq\xi^{'}$,\ for\ $+$\ is monotone,\  $\sum_{i=1}^{n}r_{b_i}^{'}\eta_{d_{i}}+\sum_{j=1}^{m}s_{c_j}^{'}\eta_{c_j}\leq\xi^{'}+\sum_{j=1}^{m}s_{c_j}^{'}\eta_{c_j}$,\ thus\ $\sum_{i=1}^{n}r_{b_i}^{'}\eta_{d_{i}}+\sum_{j=1}^{m}s_{c_j}^{'}\eta_{c_j}\leq\xi^{'}+\eta$.\ This is exactly the second condition of\ $\mathcal{D}+\eta\Rightarrow_{P}\xi +\eta$.\

To prove the continuity of\ $\cdot$,\ firstly,\ for arbitrary fixed\ $a\in\mathbb{R}^+$\ and a directed set\ $\mathcal{D}=\{\xi_i\}_{i\in I}\subseteq CX$\ with\ $\mathcal{D
}\Rightarrow_{P}\xi=\sum_{i=1}^{n}r_{b_i}\eta_{b_i}$,\ we claim that\ $a\cdot\mathcal{D}=\{a\cdot\xi_i\}_{i\in I}\Rightarrow_{P}(a\cdot\xi)$,\ since we have\ $D_i\rightarrow b_i,i=1,\dots,n$,\ and\ $\forall\ (d_1,\dots,d_n)\in\prod_{i=1}^{n},\ \forall\ r_{b_i}^{'}<r_{b_i},i=1,\dots,n$,\ there exists\ $\xi^{'}\in\mathcal{D}$,\ such that\ $\sum_{i=1}^{n}r_{b_i}^{'}\eta_{d_i}\leq\xi^{'}$\ by the hypothesis of\ $\mathcal{D}\Rightarrow_{P}\xi$.\ By the monotonicity of\ $\cdot$,\ we have\ $\sum_{i=1}^{n}a\cdot r_{b_i}^{'}\eta_{d_i}\leq a\cdot\xi^{'}$.\
Secondly,\ for arbitrary\ $\xi=\sum_{i=1}^{n}r_{b_i}\eta_{b_i}\in CX$\ and a directed set\ $D\subseteq\mathbb{R}^+$\ with\ $D\rightarrow a$.\ We claim that\ $D\cdot\xi=\{d\cdot\xi\}_{d\in D}\Rightarrow_{P}(a\cdot\xi)$.\ Let\ $D_i=\{b_i\},i=1,\dots,n$,\ $\forall\ (b_1,\dots,b_n)\in\prod_{i=1}^{n}D_i$\ and\ $\forall\ ar_{b_i}^{'}<ar_{b_i},i=1,\dots,n$,\ we only need to find\ $d\in D$\ with\ $ar_{b_i}^{'}\leq dr_{b_i},i=1,\dots,n.$\ This is accessible,\ because for each\ $i=1,\dots,n$,\ $\frac{ar_{b_i}^{'}}{r_{b_i}}<a$\ and\ $\sup D=a$.\

In conclusion,\ $(CX,\ +,\ \cdot)$\ is a directed space-cone.\ $\Box$\\

We now characterize the valuation cone as the free directed space-cone over\ $X$.\ We begin with two useful general lemmas.

\begin{proposition}\cite{GHK}\ \cite{JO1990} (Splitting Lemma )\label{spl lemma}
	For two simple valuation in\ $\mathcal{SV}(X),\ X$\ is a\ $T_0$\ space,\ we have\ $\zeta=\sum_{b\in B}r_b\eta_b\leq\sum_{c\in C}s_c\eta_c=\xi$\ if and only if there exist\ $\{t_{b,c}\in[0,+\infty):b\in B,\ c\in C\}$\ such that for each\ $b\in B,\ c\in C$,$$\sum\limits_{c\in C}t_{b,c}=r_b,\ \sum_{b\in B}t_{b,c}\leq s_c$$ and\ $t_{b,c}\neq 0$\ implies\ $b\leq c$.

\end{proposition}

\begin{lemma}\cite{GHK}\label{distinct valu}\ \cite{JO1990}
	Let\ $\zeta=\sum_{b\in B}r_b\eta_b$\ and\ $\sum_{c\in C}s_c\eta_c=\xi$\ be two simple valuations on\ $\mathcal{O}(X)$,\ where\ $X$\ is a\ $T_0$\ space.\ If\ $\xi$\ and\ $\zeta$\ are two distinct as linear combinations,\ then they are distinct as valuations.

\end{lemma}

The following theorem is the main result of this paper.

\begin{theorem}\label{theorem CX}
	Suppose\ $X$\ is a directed space,\ then\ $(CX,\ +,\ \cdot)$\ is the directed probabilistic powerspace over\ $X$,\ that is,\ endowed with topology\ $O_{\Rightarrow_P}(\mathcal{SV}(X))$,\ $(CX,\ +,\ \cdot)\cong P_P(X)$.	
\end{theorem}

{\bf Proof} Define function\ $i:X\rightarrow CX$\ as follows:\ $\forall x\in X$,\ $i(x)=\eta_x$.\ By\ \ref{distinct valu},\ $i$\ is injective.\ We prove the continuity of\ $i$.\ It is evident that if\ $x\leq y\in X$\,\ then\ $\forall\ U\in\mathcal{O}(X),\ x\in U$\ implies\ $y\in U$,\ that is\ $\eta_x\leq\eta_y$.\ Suppose we have directed set\ $D\subseteq X$\ and\ $x\in X$\ with\ $D\rightarrow x$.\ Let\ $\mathcal{D}=\{\eta_d:d\in D\}$,\ then\ $i(D)=\mathcal{D}$\ is a directed set in\ $CX$\ and\ ${\mathcal D}\Rightarrow_{P}\eta_x$.\ Since by hypothesis,\ $D\rightarrow x$,\ then\ $\forall d\in D,\ \forall k<1$,\ $k\eta_d\leq\eta_d$.\

Let\ $(Y,\ \uplus,\ \ast)$\ be a directed space-cone,\ $f:X\rightarrow Y$\ is a continuous function.\ Define\ $\bar{f}:CX\rightarrow Y$\ as follows:\ $\forall \xi=\sum_{i=1}^{n}r_{b_i}\eta_{b_i}\in CX$,\ $$\bar{f}(\xi)=\biguplus\limits_{i=1}^{n}r_{b_i}\ast f(b_i).$$
By\ Lemma\ \ref{distinct valu},\ $\bar{f}$\ is well-defined,

(1) $f=\bar{f}\circ i$.

For arbitrary\ $x\in X$,\ $(\bar{f}\circ i)(x)=\bar{f}(i(x))=\bar{f}(\eta_x)=f(x)$.

(2) $\bar{f}$\ is a directed space-cone homomorphism,\ that is,\ $\bar{f}$\ is continuous and for arbitrary\ $\sum_{i=1}^{n}r_{b_i}\eta_{b_i}$,\ $\bar{f}(\sum_{i=1}^{n}r_{b_i}\eta_{b_i})=\biguplus_{i=1}^{n}r_{b_i}\ast\bar{f}(\eta_{b_i})$.\ But this equation is evident since\ $\bar{f}(\eta_{b_i})=f(b_i),\ i=1,\dots,n.$

For continuity,\ first,\ $\bar{f}$\ is monotone.\ Let\ $\zeta=\sum_{b\in B}r_b\eta_b\leq\sum_{c\in C}s_c\eta_c=\xi$,\ by Proposition\  \ref{spl lemma}\ there exist\ $\{t_{b,c}\in[0,+\infty):b\in B,\ c\in C\}$\ such that for each\ $b\in B,\ c\in C$,$$\sum\limits_{c\in C}t_{b,c}=r_b,\ \sum_{b\in B}t_{b,c}\leq s_c$$ and\ $t_{b,c}\neq 0$\ implies\ $b\leq c$.\ From the definition of\ $\bar{f}$\ we have$$\bar{f}(\zeta)=\sum_{b\in B}r_bf(b)=\sum\limits_{b\in B}\sum\limits_{c\in C}t_{b,c}f(b)\leq\sum\limits_{c\in C}\sum\limits_{b\in B}t_{b,c}f(c)\leq\sum\limits_{c\in C}s_{c}f(c)=\bar{f}(\xi).$$

Second,\ by Proposition\ \ref{prop CXdc},\ we shall prove that\ $\bar{f}$\ perserves\ $\Rightarrow_{P}$\ convergence class, that is suppose we have a directed set\ $\mathcal{D}=\{\xi_i\}_{i\in I}\subseteq CX$\ and\ $\xi=\sum_{i=1}^{n}r_{b_i}\eta_{b_i}\in CX$\ with\ $\mathcal{D}\Rightarrow_{P}\xi$,\ then\ $\bar{f}\rightarrow\bar{f}(\xi)$\ in\ $Y$.\ $\mathcal{D}\Rightarrow_{P}\xi$\ adimits $D_i\subseteq X$\ with\ $D_i\rightarrow b_i,i=1,\dots,n$.\ For the continuity of\ $f,\ \uplus$\ and\ $\ast$,\ we have\ $$(r_{b_1}\ast f(D_1))\uplus\dots\uplus (r_{b_n}\ast f(D_n))\rightarrow\biguplus\limits_{i=1}^nr_{b_i}\ast f(b_i)=\bar{f}(\xi).$$
Here,\ $(r_{b_1}\ast f(D_1))\uplus\dots\uplus (r_{b_n}\ast f(D_n))=\{(r_{b_1}\ast f(d_1))\uplus\dots\uplus (r_{b_n}\ast f(d_n)):(d_1,\dots,d_n)\in\prod_{i=1}^{n}D_i\}$.\ For an arbitrary open neighborhood\ $U$\ of\ $\bar{f}(\xi)$,\ there exists some\ $(d_1,\dots,d_n)\in\prod_{i=1}^{n}D_i$\ such that\ $\biguplus_{i=1}^{n}r_{b_i}\ast f(d_i)\in U$.\ Fix\ $d_1,\dots,d_n$,\ and for arbitrary\ $r_{b_i}^{'}<r_{b_i},i=1,\dots,n$,\ again by the definition of\ $\mathcal{D}\Rightarrow_{P}\xi$,\ there exists some\ $\xi^{'}\in\mathcal{D}$\ such that\ $\sum_{i=1}^{n}r_{b_i}^{'}\eta_{d_i}\leq\xi^{'}$.\ Then\ $\bar{f}(\sum_{i=1}^{n}r_{b_i}^{'}\eta_{d_i})=\biguplus_{i=1}^{n}r_{b_i}^{'}f(d_i)\leq\bar{f}(\xi^{'})$.\ By the continuity of\ $f,\ \uplus$\ and\ $\ast$,\ we have\ $\biguplus_{i=1}^{n}r_{b_i}f(d_i)\leq\bar{f}(\xi^{'})$.\ But\ $U$\ is an upper set and\ $\biguplus_{i=1}^{n}r_{b_i}\ast f(d_i)\in U$,\ it follws that\ $\bar{f}(\xi^{'})\in U$,\ $\bar{f}$\ is continuous.

(3) Homomorphism\ $\bar{f}$\ is unique.

Suppose we have a directed space-cone homomorphism\ $g:(CX,\ +,\ \cdot)\rightarrow (Y,\ \uplus,\ \ast)$\ such that\ $f=g\circ i$,\ then\ $g(\eta_x)=f(x)=\bar{f}(\eta_x)$.\ For each\ $ \xi=\sum_{i=1}^{n}r_{b_i}\eta_{b_i}\in CX$,\ \begin{eqnarray*}
 g(\xi)&=& g(r_{b_1}\eta_{b_1}+r_{b_2}\eta_{b_2}\cdots+r_{b_n}\eta_{b_n})\\ &=&g(r_{b_1}\eta_{b_1})\uplus g(r_{b_2}\eta_{b_2})\uplus \cdots \uplus g(r_{b_n}\eta_{b_n})\\ &=&r_{b_1}\ast g(\eta_{b_1})\uplus r_{b_2}\ast g(\eta_{b_2})\uplus \cdots\uplus r_{b_n}\ast g( \eta_{b_n})\\ &=&r_{b_1}\ast f({b_1})\uplus r_{b_2}\ast f({b_2})\uplus \cdots\uplus r_{b_n}\ast f({b_n})\\ &=&\bar{f}(\xi).
 \end{eqnarray*}
 Thus\ $\bar{f}$\ is unique.

 In conclusion,\  according to definition\ \ref{def dpp},\ endowed with topology\ $O_{\Rightarrow_P}(\mathcal{SV}(X))$,\ the directed space-cone\ $(CX,\ +,\ \cdot)$\ is the directed probabilistic powerspace of\ $X$,\ that is,\ $P_P(X)\cong (CX,\ +,\ \cdot)$.\ $\Box$

 \vskip 3mm
The directed probabilistic powerspace is unique in the sense of algebraic isomorphism and topological homeomorphism,\  so we can directly denote the directed probabilistic powerspace by\ $P_P(X)=(CX,\ +,\ \cdot)$\ of each directed space\ $X$.

 Suppose\ $X,\ Y$\ are two directed spaces,\ $f:X\rightarrow Y$\ is a continuous function.\ Define function\ $P_P(f):P_P(X)\rightarrow P_P(Y)$\ as follows:\ $\forall \xi=\sum_{i=1}^{n}r_{b_i}\eta_{b_i}\in CX$,\ $$P_P(f)(\xi)=\sum_{i=1}^{n}r_{b_i}\eta_{f(b_i)}.$$
 We can check that,\ $P_P(f)$\ is well-defined and order preserving.\ It is easy to check that,\ $P_P(f)$\ is a directed space-cone homomorphism between these two extended probabilistic powerspaces.\ If\ $id_X$\ is the identity function
and\ $g:Y\rightarrow Z$\ is an arbitrary continuous function from\ $Y$\ to a directed space\ $Z$,\ then,\ $P_P(id_X)=id_{P_P(X)},\ P_P(g\circ f)=P_P(g)\circ P_P(f)$.\ Thus,\ $P_P:{\bf Dtop}\rightarrow {\bf Dcone}$\ is a functor from \ ${\bf Dtop}$\ to\ ${\bf Dcone}$.\ Let\ $U:{\bf Dcone}\rightarrow {\bf Dtop}$\ be the forgetful functor. By theorem\ \ref{theorem CX},\ we have the following result.

 \begin{corollary}
 $P_P$\ is a left adjoint of the forgetful functor\ $U$,\ that is,\ ${\bf Dcone}$\ is a reflective subcategory of\ ${\bf Dtop}$.
 \end{corollary}

\section{Relations Between Directed Probabilistic Powerspace And Probabilistic Powerdomain}

In this section,\ we will discuss the relations between the extended probabilistic powerdomain of dcpo and the directed probabilistic powerspace.

According to the results in the last section,\ for an arbitrary directed space\ $X$,\ the directed probabilistic powerspace is the set of\ $CX$\ endowed with the\ $\Rightarrow_P$\ convergence topology.\ In general,\ for an arbitrary topological space and arbitrary\ dcpo,\ although the extended probabilistic powerdomain exists,\ their concrete structure cannot be expressed explicitly(see article\ \cite{BS2015, Hec91, Hec92}).

\begin{definition}{\rm\cite{GHK}}
	For a topological space\ $X$,\ the valuation powerdomain(also called the extended probabilistic powerdomain)\ $\mathcal{V}(X)$\ of\ $X$\ is the set of all continuous valuations on\ $\mathcal{O}(X)$\ with the pointwise order.

\end{definition}

\begin{definition}{\rm\cite{GHK}\ \cite{JO1990}}\label{def prec}
	For a simple valuation\ $\xi=\sum_{b\in B}r_b\eta_b$\ and a continuous valuation\ $\mu$\ on\ $X$,\ a domain equipped with the Scott topology,\ we set\ $\xi\prec\mu$\ if for all nonempty\ $K\subseteq B$,\ we have\ $\sum_{b\in K}r_b<\mu(\bigcup_{b\in K}\dua b).$

\end{definition}

\begin{theorem}{\rm\cite{GHK}\ \cite{JO1990}}\label{thm prec}
	For a domain\ $X$\ the valuation powerdomain\ $\mathcal{V}(X)$\ is a domain.\ Each continuous valuation\ $\mu$\ is the directed supremum of the simple valuations way below it,\ and for simple valuation\ $\xi$,\ one has\ $\xi\ll\mu\iff\xi\prec\mu$.

\end{theorem}

\begin{theorem}{\rm\cite{GHK}\ \cite{JO1990}}
	Given any dcpo-cone\ $C$\ and a continuous function\ $f:X\rightarrow C$,\ where\ $X$\ is a domain equipped with the Scott topology,\ there exists a unique continuous linear map\ $f^{\ast}:\mathcal{V}(X)\rightarrow C$\ such that\ $f^{\ast}\eta_X=f$,\ here\ $\eta_X(x)=\eta_x,\ \forall x\in X$.

\end{theorem}

\begin{definition}\label{def llcurly}
	Suppose\ $X$\ is a topological space,\ $\mu,\ \nu\in\mathcal{SV}(X)$,\ and\ $\mu=\sum_{b\in B}r_b\eta_b,\ \nu=\sum_{c\in C}s_c\eta_c$,\ we say\ $\mu\llcurly\nu$\ if there exist\ $\{t_{b,c}\in[0,+\infty):b\in B,\ c\in C\}$\ such that for each\ $b\in B,\ c\in C$,$$\sum\limits_{c\in C}t_{b,c}=r_b,\ \sum_{b\in B}t_{b,c}< s_c$$ and\ $t_{b,c}\neq 0$\ implies\ $c\in\interior(\uparrow b)$.

\end{definition}

Since the soberification of each c-space is a continuous domain,\ by Theorem \uppercase\expandafter{\romannumeral4}-9.16, Proposition\ \uppercase\expandafter{\romannumeral4}-9.18,\ and\ \uppercase\expandafter{\romannumeral4}-9.19\ in\cite{GHK},\ we have the following propositions.

\begin{proposition}\label{prop llcurly}
	Suppose\ $X$\ is a c-space,\ $\mu,\ \nu,\ \xi\in\mathcal{SV}(X)$,\ then the followings hold:
	\begin{enumerate}
		\item $\mu\llcurly\nu\ \Rightarrow\ \mu\leq\nu$.
		\item $\xi\llcurly\mu\leq\nu\ \Rightarrow \xi\llcurly\nu$.
		\item $\mu,\ \nu\llcurly\xi\ \Rightarrow\exists\xi^{'}\in\mathcal{SV}(X),\ \mu,\ \nu \llcurly \xi^{'} \llcurly \xi$.
		\item $\mu\nleq\nu\ \Rightarrow\exists\xi^{'}\in\mathcal{SV},\ \xi^{'}\llcurly\mu\ \&\ \xi^{'}\nleq\nu$.
	\end{enumerate}
\end{proposition}

 According to the definition\ \ref{def llcurly},\ if we have\ $\mu\llcurly\nu$,\ then each\ $b\in B$,\ there exists some\ $c\in C$\ such that\ $c\in\interior(\uparrow b)$.\ Let\ $\Uparrow\mu=\{\nu\in\mathcal{SV}(X):\mu\llcurly\nu\}$.\ We claim that\ $\Uparrow\mu$\ is an open set in\ $CX$,\ and\ $CX$\ is a c-space.

 \begin{lemma}\label{lemma llopen}
 	Let\ $X$\ be a c-space,\ then for each\ $\mu\in\mathcal{SV}(X)$,\ $\Uparrow\mu$\ is open in\ $CX$.
 	
 \end{lemma}
\noindent{\bf Proof}  We only need to check that\ $\Uparrow\mu$\ is a\ $\Rightarrow_P$\ convergence open set.\ Suppose we have a directed set\ $\mathcal{D}\subseteq\mathcal{SV}(X)$\ and\ $\xi\in\mathcal{SV}(X)$\ with\ $\mathcal{D}\Rightarrow_p\xi\in\Uparrow\mu$.\ Let\ $\mu=\sum_{i=1}^{k}r_{b_i}\eta_{b_i}$,\ with support\ $B$,\ $\xi=\sum_{j=1}^{n}s_{c_j}\eta_{c_j}$,\ with support\ $C$.\ By the definition of\ $\Rightarrow_{P}$\ convergence,\ there exist directed sets\ $D_c\subseteq X$\ such that\ $D_c\rightarrow c,\ c\in C.$\ Suppose we have\ $C\cap(\bigcup_{b\in B}\interior(\uparrow b))=\{c_1,\cdots,c_k\}=K,\ 1\leq k\leq n.$\ For each\ $c\in K$\, we have a finite set\ $B_c\subseteq B$\ such that\ $\forall b\in B_c,\ c\in\interior(\uparrow b)$.\ Since\ $D_c\rightarrow c$,\ and  eventually in\ $\interior(\uparrow b)$,\ we have finitely\ $d_{cb}\in\interior(\uparrow b),\ b\in B_c$.\ Then we may pick the largest one\ $d_c\in\bigcap_{b\in B_c}\interior(\uparrow b)$.\ By hypothesis,\ there exist\ $\{t_{b,c}\in[0,+\infty):b\in B,\ c\in C\}$\ satisfy the definition of\ $\mu\llcurly\xi$.\ Define a new valuation\ $\nu=\sum_{c\in C}s_c\eta_{d_c}$.\ It is direcly to check that\ $\mu\llcurly\nu$\ with $\{t_{b,c}\in[0,+\infty):b\in B,\ c\in C\}$.\ Again by the definition of\ $\mathcal{D}\Rightarrow_{P}\xi$,\ there exists some\ $\xi^{'}\in\mathcal{D}$\ such that\ $\nu\leq\xi^{'}$.\ Now we have\ $\mu\llcurly\nu\leq\xi^{'}$,\ according to\ 2\ in Proposition\ \ref{prop llcurly},\ $\mu\llcurly\xi^{'}$.\ It follows that\ $\xi^{'}\in\mathcal{D}\cap\Uparrow\mu$,\ $\Uparrow\mu$\ is open.\ $\Box$\\

\begin{theorem}
	Let\ $X$\ be a c-space,\ then\ $CX$\ is a c-space.
	
\end{theorem}
\noindent{\bf Proof} Suppose we have\ $\xi=\sum_{i=1}^{n}s_{c_i}\eta_{c_i}\in\mathcal{U}$,\ where\ $\mathcal{U}$\ is a\ $\Rightarrow_{P}$\ convergence open set.\ Since\ $X$\ is a c-space,\ then\ $\mathcal{D}=\{\sum_{i=1}^{n}r_i\eta_{b_i}:r_i<s_{c_i},\ c_i\in\interior(\uparrow b_i)\}$\ is directed and\ $\mathcal{D}\Rightarrow_P\xi$.\ Thus we have some\ $\xi^{'}\in\mathcal{D}\cap\mathcal{U}$,\ and it is directly to check that each\ $\mu\in\mathcal{D}$,\ $\mu\llcurly\xi$.\ Since\ $\{b_{b_{i}c_{j}}:t_{b_{i}c_{j}}=0,i\neq j;t_{b_{i}c_{j}}=1,i=j\}$\ satisfy the definition of\ $\mu\llcurly\xi$.\ By lemma\ \ref{lemma llopen},\ $\Uparrow\xi^{'}$\ is open and\ $\xi\in\Uparrow\xi^{'}\subseteq\mathcal{U}$,\ $CX$\ is a c-space.\ $\Box$\\



Limited to the space,\ we may not generalize all the auxiliary relations from domains to c-spaces.\ Thus,\ we have the following proposition limited to domains.\ LetA\ $X$\ be a continuous\ domain.\ Then\ $(X,\ \sigma(X))$\ is a directed space,\ and\ $CX\subseteq \mathcal{V}(X)$.\ By\ $\sigma(\mathcal{V}(X))|_{\mathcal{SV}(X)}$\ denote the relative topology from the Scott topology on\ $\mathcal{SV}(X)$.

\begin{proposition}\label{prop rstCX}
	Suppose\ $X$\ is a domain,\ then $O_{\Rightarrow_P}(\mathcal{SV}(X))=\sigma(\mathcal{V}(X))|_{\mathcal{SV}(X)}$,\ that is\ $CX$\ is the subspace of\ $(\mathcal{V}(X),\ \sigma(\mathcal{V}(X))$.
	
\end{proposition}

\noindent{\bf Proof}
\noindent  Suppose\ $\mathcal{U}\in O_{\Rightarrow_P}(\mathcal{SV}(X))$.\ Let\ $\mathcal{U}_V=\{\mu\in \mathcal{V}(X):\exists\xi\in {\mathcal U},\ \xi\leq \mu\}$.\ Obviously,\ ${\mathcal U}= \mathcal {U}_V\cap CX$.\ We claim that\ $\mathcal {U}_V$\ is Scott open.\ Suppose we have a directed set of continuous valuations\ $\mathcal F=\{\mu_i\}_{i\in I}\subseteq\mathcal{V}(X)$\ with\ $\sup_{i\in I}\{\mu_i\}=\mu\in\mathcal{U}_V$.\ By the definition of\ $\mathcal{U}_V$,\ there exists some\ $\xi=\sum_{i=1}^{n}r_{b_i}\eta_{b_i}\in\mathcal{U}$\ such that\ $\xi\leq\mu$.\ Since\ $X$\ is a continuous domain,\ we have finitely directed sets\ $D_i=\dda b_i$\ with\ $D_i\rightarrow b_i,i=1,\dots,n$.\ It is straightly to verify that\ $\mathcal D=\{\sum_{i=1}^{n}r_{i}\eta_{a_i}:r_i<r_{b_i},a_i\ll b_i,i=1,\dots,n\}$\ is directed and\ $\mathcal D\Rightarrow_P\xi$.\ By the hypothesis of\ $\mathcal U$,\ there exists some\ $\xi^{'}\in\mathcal D\cap\mathcal U$.\ According to the definition of\ $\mathcal{D}$\ and definiton\ \ref{def prec},\ we have\ $\xi^{'}\prec\xi$,\ then,\ by \ref{thm prec},\ $\xi^{'}\ll\xi$.\ $\mathcal{V}(X)$\ is a continuous domain implies\ $\dua\xi^{'}$\ is a Scott open set.\
Since\ $\mathcal F=\{\mu_i\}_{i\in I}$\ is a directed set with supremum\ $\mu\in\dua \xi^{'}$,\ then we have some\ $\mu^{'}\in\mathcal F$\ with\ $\xi^{'}\ll\mu^{'}$,\ thus\ $\xi^{'}\leq\mu^{'}$\ and\ $\mu^{'}\in\mathcal{U}_V$,\ it follows that\ $\mathcal {U}_V$\ is Scott open.\

On the other hand,\ let\ ${\mathcal V}\in\sigma(\mathcal{V}(X))$,\ and\ ${\mathcal D}\subseteq CX$\ is a directed set with\ ${\mathcal D}\Rightarrow_P\xi=\sum_{i=1}^{n}r_{b_i}\eta_{b_i}\in {\mathcal V}\cap CX$.\ We claim that\ $\mathcal{D}\Rightarrow_{P}\xi$\ implies\ $\xi\leq\sup\mathcal{D}$.\ By the definition,\ we have finitely directed sets\ $D_i\subseteq X$\ with\ $D_i\rightarrow b_i,i=1,\dots,n$.\ For arbitrary open set\ $U\in\mathcal{O}(X)$\ and suppose\ $b_1,\dots,b_k\in U,0\leq k\leq n$. Thus we have\ $(d_1,\dots,d_n)\in\prod_{i=1}^{k}D_i$\ such that\ $d_i\in U,i=1,\dots,k$.\ Again by the definition of\ $\mathcal{D}\Rightarrow_P\xi$,\ $\forall r_{b_i}^{'}<r_{b_i},i=1,\dots,k$,\ there exists some\ $\xi^{'}\in\mathcal{D}$\ with\ $\sum_{i=1}^{k}r_{b_i}^{'}\eta_{d_i}\leq\xi^{'}$,\ then\ $$\xi(U)=(\sum_{i=1}^{k}r_{b_i}^{'}\eta_{d_i})(U)\leq\xi^{'}(U)\leq (\sup\mathcal{D})(U).$$
It follows that\ $\sup\mathcal{D}\in\mathcal{V}$,\ thus there exists some\ $\xi^{''}\in\mathcal{D}\cap\mathcal{V} $,\ then\ $\xi^{''}\in\mathcal{D}\cap\mathcal{V}\cap CX$,\ $\mathcal{V}\cap CX$\ is open in\ $CX$.\ $\Box$\\

Finally,\ we ending this paper with an example.
\begin{example}\label{example CXneq}
	
	Let\ $X=[0,1]$\ endowed with the Scott topology,\ $\forall a\in [0,1],(a,1]\in\mathcal{O}(X)$,\ let\ $\mu:\mathcal{O}(X)\rightarrow \mathbb{R}^{+},\mu((a,1])=1-a$.\ Then\ $\mu$\ is a Scott continuous valuation but\ $\mu\notin\mathcal{SV}(X)$,\ since each\ $\xi=\sum_{i=1}^{n}r_{b_i}\eta_{b_i}\in\mathcal{SV}(X)$,\ the range of\ $\eta_{b_i},i=1,\dots,n,$\ is\ $0$\ or\ $1$,\ it follows that the range of\ $\xi$\ is finite.\ But the range of\ $\mu$\ is\ $[0,1]$.\ Thus the extended probabilistic powerspace of\ $X$\ is not equal to its extended probabilistic powerdomain.
	
\end{example}

\section*{Acknowledgement}
We thank the anonymous referee for improving the presentation of the paper,\ we would like to thank Chen Yu and Chen Yuxu for patiently and helpful discussions.


\section*{Reference}

\begin{thebibliography}{99}
\bibitem{AJ} Abramsky, S.,Jung, A.:\ Domain theory. In: Abramsky, S.,Gabbay, D.M., Maibaum, T.S.E.(eds.),  Semantic Structures. In: Handbook of Logic in Computer Science,  vol.3, pp.1-168, Clarendon Press, Oxford (1994)


\bibitem{BS2015} Battenfeld, I., Sch\"{o}der, M.:Observationally-induced lower and upper powerspace constructions. Journal of Logical and Algebraic Methods in Programming. 84, 668-682 (2015)
\bibitem{CS1991} Erne E. The ABC of order and topology. In: H. Herrlich, H. E. Porst (des.) Category Theory at Work, pp. 57-83, Heldermann, Berlin (1991)

\bibitem{EN2009}Ern\'{e}, M.: Infinite distributive laws versus local connectedness and compactness properties. Topology and its Applications. 156, 2054-2069 (2009)

\bibitem{GK2017}Geng J., Kou H.: Consistent Hoare powerdomains over dcpos,  Topology and its Applications. 232, 169-175 (2017)

\bibitem{GHK}Gierz, G. et al.:  Continuous Lattices and Domains.  Cambridge University Press,  Cambridge (2003).

\bibitem{GouJ14} Goubault-Larrecq, J., Jung, A.: QRB, QFS, and the Probabilistic Powerdomain. Electronic Notes in Theoretical Computer Science. 308, 167-182 (2014)

\bibitem{Gou13}Goubault-Larrecq, J.:  Non-Hausdorff topology and domain theory: Selected topics in point-set topology. Cambridge University Press, Cambridge (2013)

\bibitem{GRA1987}Graham, S.K.: Closure properties of a probabilistic domain construction. Springer, Berlin, Heidelberg. 213-233 (1987)

\bibitem{HACK1993} Heckmann, R.: Probabilistic power domains,information systems, and locales. Springer, Berlin, Heidelberg. 410-437 (1993)

\bibitem{HAC1995} Heckmann, R.: Spaces of valuations. In Andima, S., Flagg, R.C., Itzkowitz, G., Misra,  P., Kong,  Y. and
Kopperman,R., editors, Papers on General Topology and Applications: Eleventh Summer Conference at the University of Southern
Maine, volume 806 of Annals of the New York Academy of Sciences. 174-200 (1996)

\bibitem{Hec91}Heckmann, R.: Power domain\ constructions. Sci. Comput. Program. 17, 77-117 (1991)

\bibitem{Hec92}Heckmann, R.:  An upper power domain\ construction in terms of strongly compact sets,  Lect. Notes Comput. Sci. 598, 272-293 (1992)

\bibitem{HK13}Heckmann, R., Keimel, K.:  Quasicontinuous domains and the Smyth powerdomain,  Electron. Notes Theor. Comput. Sci. 298, 215-232 (2013)

\bibitem{PLOT1989}Jones, C.G., Plotkin, D.: A probabilistic powerdomain of evaluations. Fourth Annual Symposium on Logic in Computer Science. IEEE Computer Society. 186-195 (1989)

\bibitem{JO1990} Jones, C.: Probabilistic Non-Determinism, Ph.D. Thesis, University of Edinburgh, Report ECS-LFCS-90-105(1990)

\bibitem{JT98}Jung, A., Tix, R.: The troublesome probabilistic powerdomain. Electronic Notes in Theoretical Computer Science. 13, 70-91 (1998)

\bibitem{KEI2006}Keimel,K.:\ Topological cones: functional analysis in a $T_0$-setting. Springer-Verlag. 77(1), 109-142 (2008)

 \bibitem{Kou15}Kou, H.: Directed spaces: An extended framework for domain\ theory, 1th Pan Pacific International Conference on Topology and Applications£¬Min Nan Normal University,  Zhangzhou City(2015). 11. 25-30

 \bibitem{LyuK2018}Lyu Z., Kou H.: The probabilistic powerdomain from a topological viewpoint. Topology and its Applications. 237, 26-36. (2018)

 \bibitem{Mac71}MacLane, S.: Categories for the Working Mathematician. Springer-Verlag, New York (1971)

\bibitem{GEN}  Mislove,M.:\ Generalizing domain theory.\ International Conference on Foundations of Software Science and Computation Structure. Springer, Berlin, Heidelberg. 1-19 (1998)


\bibitem{TND}Mislove, M.:\ Topology,domain theory and theoretical computer science.Topology and its Applications. 89(1-2), 3-59 (1998)

\bibitem{Mis88}Mislove, M.: On the Smyth power domain. Lect. Notes Comput. Sci. 298, 161-172 (1988)

\bibitem{Plo76}Plotkin, G.D.: A powerdomain construction.\ SIAM J. Comput. 5, 452-487 (1976)

\bibitem{SD1980} Saheb-Djahromi,N.:CPO's of measures for nondeterminism. Theoretical Computer Science. 12(1), 19-37 (1980)

\bibitem{Scott70}Scott, D. S.: Outline of a mathematical theory of computation. In 4th Annual Princeton Conference on Information Sciences and Systems.(1970)

\bibitem{Scott71}Scott, D. S.:  Continuous Lattices: Toposes, Algebraic Geometry and Logic. Springer Lecture Notes in Mathematics. 274, 97-136 (1972)

\bibitem{Scott82}Scott, D. S.\ Lectures on a mathematical theory of computation,\ In M.Broy and G.Schmidt,\ editors,\ Theoretical Foundations of Programming Methodology. 145-292 (1982)

\bibitem{SMY}Smyth, M. B.: Power domains and predicate transformers: A topological view. Springer, Berlin, Heidelberg. 662-675 (1983)

\bibitem{TKP2009}Tix, R.,Keimel, K., Plotkin, G.: Semantic domains for combining probability and non-determinism. Electronic Notes in Theoretical Computer Science. 222, 3-99(2009)

\bibitem{XXLK2020}Xie, X., Kou H.:\ Lower power structures of directed spaces(Chinese).\ Journal of Sichuan University (Natural Science Edition). 57, 211-217 (2020)

\bibitem{Xu16} Xu X.:\ Order and Topology(Chinese).\ Science Press, Beijing\ (2016)

\bibitem{YYK2014} Yu, Y.,Kou, H.:\ On directed Space(Chinese).\ Journal of Sichuan University (Natural Science Edition) (2014)

\bibitem{YYK2015} Yu Y., Kou, H.:\ Directed spaces defined through\ $T_0$\ spaces with specialization order(Chinese),\ Journal of Sichuan University (Natural Science Edition). 52(2), 217-222 (2015)

\bibitem{YK2014-1}Yuan, Y., Kou, H.:\ Consistent Smyth Powerdomains.\ Topology and its Applications. 173, 264-275 (2014)

\bibitem{YK2014-2} Yuan, Y., Kou, H.:\ Consistent Hoare Powerdomains.\ Topology and its Applications.\ 178, 40-45 ( 2014)

\bibitem{YK2014-3} Yuan, Y., Kou, H.:\ Consistent Plotkin powerdomains.\ Topology and its Applications.\ 178, 339-344 (2014)






\end{thebibliography}
\end{document}